\documentclass[12pt]{article}

\usepackage{CJK,CJKnumb,amsmath}
\usepackage{amsfonts}
\usepackage{amsthm}
\usepackage{geometry}
\usepackage{mathrsfs}
\usepackage{amsmath}
\usepackage{amssymb}
\usepackage{amsfonts}
\usepackage[percent]{overpic}
\usepackage{bm}
\usepackage[all]{xy}
\usepackage{graphicx}
\usepackage{subfigure}
\usepackage{latexsym}
\usepackage[colorlinks, linkcolor=blue, anchorcolor=bule, citecolor=blue]{hyperref}
\usepackage{epigraph}
\usepackage{hyperref}
\usepackage{fancyhdr}
\usepackage{comment}
\usepackage[utf8]{inputenc}

\usepackage{mathtools}
\mathtoolsset{showonlyrefs}

\numberwithin{equation}{section}
\usepackage{color}

\newtheorem{theorem}{Theorem}

\theoremstyle{remark}
\newtheorem{remark}{Remark}

\newtheorem*{ack}{Acknowledgement}

\def\area{\operatorname{area}}
\def\mod{\operatorname{mod}}

\def\sing{\operatorname{Sing}}
\def\dens{\operatorname{dens}}

\def\diam{\operatorname{diam}}

\def\dim{\operatorname{dim}}
\def\hc{\operatorname{\widehat{\mathbb C}}}
\def\c{\operatorname{\mathbb C}}
\def\N{\operatorname{\mathbb N}}

\def\dim{\operatorname{dim}}

\def\j{\operatorname{\mathcal{J}}}
\def\f{\operatorname{\mathcal{F}}}

\def\I{\operatorname{\mathcal{I}}}
\def\b{\operatorname{\mathcal{B}}}
\def\s{\operatorname{\mathcal{S}}}

\def\sing{\operatorname{Sing}}
\def\res{\operatorname{Res}}

\begin{document}
\title{Hausdorff dimension of escaping sets of Nevanlinna functions}
\author{Weiwei Cui}
\date{}

\maketitle

\begin{center}
\emph{Dedicated to Professor Dr. Walter Bergweiler on the occasion of his 60th birthday}
\end{center}

\begin{abstract}
We determine the exact values of Hausdorff dimensions of escaping sets of meromorphic functions with polynomial Schwarzian derivatives. This will follow from the relation between these functions and the second order differential equations in the complex plane.

\medskip
\emph{Mathematics Subject Classification}: 30D05 (primary), 37F10 (secondary).

\medskip
\emph{Keywords}: Meromorphic functions, Schwarzian derivatives, Hausdorff dimension, escaping sets.

\end{abstract}

\section{Introduction and main result}

Let $f:\c\to\hc$ be a transcendental meromorphic function. The \emph{Fatou set} $\f(f)$ of $f$ is defined as the set of points in $\c$ where the iterates are defined and form a normal family. The \emph{Julia set} $\j(f)$ is the complement of $\f(f)$ in $\c$. These two sets are the main objects in transcendental dynamics. We refer to \cite{bergweiler1} for an introduction to transcendental dynamics.

Another concept, which has attracted a lot of interest recently, is the \emph{escaping set} $\I(f)$, defined as
$$\I(f):=\left\{z\in\c:\, f^{n}(z)\to\infty~\text{as}~n\to\infty\right\}.$$
This set was introduced into transcendental dynamics by Eremenko \cite{eremenko3}. McMullen \cite{mcmullen11} proved that  the two-dimensional Lebesgue measure of $\I(\sin(\alpha z+\beta))$ is positive and that the Hausdorff dimension of $\I(\lambda e^z)$ is two, where $\alpha,\beta, \lambda\in\c$ and $\alpha, \lambda\neq 0$. For further results concerning the Lebesgue measure and Hausdorff dimension of Julia sets and escaping sets of entire functions we refer to \cite{baranski1, schubert, bergweiler10, eremenko2, aspenberg1, cuiwei1} and the references therein. 

For meromorphic functions with a direct singularity over infinity the situation is similar to entire functions; cf. \cite{bergweiler3}. However, for meromorphic functions for which $\infty$ is not an asymptotic value, the situation is different; see \cite[Theorem 1.1]{bergweiler2} for details. For further results concerning the Hausdorff dimensions of Julia sets and escaping sets, we refer to \cite{stallard3}.

The main purpose of this paper is to compute the Hausdorff dimension of the escaping set of meromorphic functions with polynomial Schwarzian derivative. Recall here that the \emph{Schwarzian derivative} of $f$ is defined by
$$\s_{f}=\left(\frac{f''}{f'} \right)'-\frac{1}{2}\left(\frac{f''}{f'} \right)^2.$$
Meromorphic functions whose Schwarzian derivatives are entire functions have close connections with the second order differential equations in complex domains; see \cite[Chapter 5]{hille3} and also the following section. In particular, these meromorphic functions have the property that the spherical derivative of $f$ is never zero. In other words, they are locally univalent.

In this paper, we are mainly interested in meromorphic functions whose Schwarzian derivatives are polynomials. Nevanlinna was the first to explore these functions in terms of the inverse problem in the value distribution theory of meromorphic functions \cite{nevanlinna6, goldbergmero}. Following \cite{eremenko6} we say that a  meromorphic function is called a \emph{Nevanlinna function} if it has a polynomial Schwarzian derivative. Every such function has only finitely many asymptotic values. Many familiar functions are Nevanlinna functions, including, for instance, the tangent function $\lambda \tan(z)$, the exponential function $\lambda e^z$ and also functions of the form $\int^{z}_{0}\exp({p(s)})ds$ with a polynomial $p$.

\medskip
Nevanlinna functions appeared in the setting of transcendental dynamics first in \cite{devaney16}, in which Devaney and Keen showed that these functions have no wandering domains, which is now known for all meromorphic functions with finitely many singular values \cite{baker6}. Qiu proved in \cite{qiu1} that the Hausdorff dimension of the Julia set and the escaping set of a Nevanlinna function for which $\infty$ is an asymptotic value is two. This has been generalized to meromorphic functions of finite order with a bounded set of singular values for which infinity is an asymptotic value \cite{baranski1, schubert, bergweiler3}.  Bergweiler and Kotus \cite{bergweiler2} proved that if a meromorphic function has a bounded set of singular values and finite order, and infinity is not an asymptotic value, then the Hausdorff dimension of the escaping set is strictly less than $2$. In fact, they gave an upper bound in terms of the order of growth and the maximal multiplicity of the poles.

\medskip
Our theorem below shows that in case that $\infty$ is not an asymptotic value of Nevanlinna functions, the Hausdorff dimension of the escaping set is strictly less than one and can in fact be determined explicitly. This value, as one can see in the following, depends only on the degree of the corresponding Schwarzian derivatives rather than any dynamical assumptions.

\medskip
To prove our result, we will consider the following set, for $R>0$,
$$\j_{R}(f):=\left\{z\in\c:~|f^{n}(z)|>R~\text{for all}~n\in\mathbb{N}\right\}.$$
Note that every escaping point will finally fall into this set. The Hausdorff dimension of a set $A$ will be denoted by $\dim A$. Now we state our theorem as follows.

\begin{theorem}\label{nq}
Let $f$ be a Nevanlinna function for which $\infty$ is not an asymptotic value. Then 
$$\dim\I(f)=\lim_{R\to\infty}\dim\j_{R}(f)=\frac{m+2}{m+4},$$
where $m$ denotes the degree of the Schwarzian derivative of $f$.
\end{theorem}

Combined with previous discussions, our result gives a complete classification of Nevanlinna functions in terms of their Hausdorff dimensions of escaping sets: For Nevanlinna functions, the Hausdorff dimension of the escaping set is either $2$ or $(m+2)/(m+4)$ with $m$ as above.

\smallskip
In the context considered here, it follows from \cite[Theorem 1.1]{bergweiler2} that $\dim\I(f)\leq (2m+4)/(m+6)$. Thus our result sharpens their result for Nevanlinna functions. Moreover, the theorem also generalizes known results on certain explicit functions. For instance, it is known before that the $\dim\I(f)=1/2$ if $f(z)=\lambda\tan(z)$ for $\lambda\neq 0$, see \cite{kotus2} for instance.

\begin{remark}
A result of Mayer implies that the Hausdorff dimension of the Julia set of a Nevanlinna function is strictly larger than that of the escaping set \cite{mayer2}. Mayer actually proved that for Nevanlinna functions, the hyperbolic dimension of the function is strictly larger than ${(m+2)}/{(m+4)}$, which, by the above theorem, is equal to the Hausdorff dimension of the escaping sets. Examples of transcendental \emph{entire} functions with the property that the Hausdorff dimension of the Julia set is strictly larger than that of escaping set were constructed by Rempe-Gillen and Stallard \cite[Theorem 1.1]{rempe11}.
\end{remark}

\begin{remark}
The topological completeness of Nevanlinna functions implies that any two topologically equivalent Nevanlinna functions in the sense of \cite{eremenko2} have escaping sets of the same Hausdorff dimension. It is an open problem whether there exist quasiconformally equivalent \emph{meromorphic} functions in the class $\b$ whose Hausdorff dimensions of the escaping sets are differrent \cite{bergweiler10, rempe11}. 
\end{remark}

The result may be extended to meromorphic functions whose Schwarzian derivatives are rational functions. Such functions have only finitely many critical points and finitely many asymptotic values. Instead of considering a second order complex differential equation with a polynomial coefficient, one need in this case consider rational function coefficients. Our proof goes through in this case, but we leave the details to the reader. For future references, we state this result as a theorem.

\begin{theorem}
Let $f$ be a meromorphic function whose Schwarzian derivative is a rational function of degree $m$. If $\infty$ is not an asymptotic value of $f$, then
$$\dim\I(f)=\lim_{R\to\infty}\dim\j_{R}(f)=\frac{m+2}{m+4}.$$
\end{theorem}

\bigskip
A meromorphic function belongs to the class $\b$ if the set of finite singular values is bounded. See the next section for precise definitions. The proof of the above theorems works in a more general setting and actually yields the following result.

\begin{theorem}
Let $f\in\b$ be a meromorphic function for which $\infty$ is not a singular value. Suppose that $\{a_j\}$ are poles of $f$ ordered in such a way that $\dots\leq|a_j|\leq|a_{j+1}|\leq\dots$. Suppose furthermore that $a_j$ and its residue $b_j$ satisfy
$$|a_j|\sim j^{2/(m+2)}\,\,\,\,\text{and}\,\,\,~|b_j|\sim j^{-m/(m+2)}$$
for $j$ large. Then $$\dim\I(f)=\lim_{R\to\infty}\dim\j_{R}(f)=\frac{m+2}{m+4}.$$
\end{theorem}

\medskip
\noindent{\emph{Notations}.} The closure of a set $A$ in the plane is denoted by $\overline{A}$. The Euclidean diameter and spherical diameter of a set $A$ will be denoted by $\diam A$ and $\diam_{\chi} A$, respectively. The residue of $f$ at a pole $z$ is written as $\res(f,z)$.  A Euclidean disk centred at $z$ of radius $r$ will be denoted by $D(z,r)$ while $D_{\chi}(z,r)$ means spherical disk. The Euclidean (resp. spherical) area of a set $A$ on the plane is denoted by $\area(A)$ (resp. $\area_{\chi}(A)$). Let $A,\,B$ be measurable sets of the plane. Then the density and the spherical density of $A$ in $B$, denoted respectively by $\dens(A,B)$ and $\dens_{\chi}(A,B)$, are defined as 
$$\dens(A, B)=\frac{\area (A\cap B)}{\area (B)},~~~\,\,\,~~~~\dens_{\chi}(A,B)=\frac{\area_{\chi}(A\cap B)}{\area_{\chi}(B)}.$$

\section{Preliminaries}

The set $\sing(f^{-1})$ of singularities of the inverse of a meromorphic function $f$ is defined as the set of all finite critical and asymptotic values of $f$. (For transcendental entire functions, $\infty$ is always an asymptotic value, while for meromorphic functions, $\infty$ can serve as an asymptotic, critical or regular value.) The \emph{Eremenko-Lyubich class} $\b$ consists of meromorphic functions with a bounded singular set $\sing(f^{-1})$. If, in particular, $\sing(f^{-1})$ is finite, then we say that $f$ belongs to the \emph{Speiser class} $\s$. These classes of functions were considered in dynamics first by Eremenko and Lyubich \cite{eremenko2} and have attracted much interest since then, mainly due to the fact that these functions have certain dynamical properties similar to those of polynomials.

Let $f$ be a meromorphic function in the plane. It is known that $\I(f)\neq\emptyset$ and $\j(f)=\partial\I(f)$. The transcendental entire case was proved by Eremenko \cite{eremenko3} and the general meromorphic case is due to \cite{dominguez1}. If, in addition, $f\in\b$, then $\I(f)\subset\j(f)$ \cite{eremenko2, rippon6}. Therefore, we have $\j(f)=\overline{\I(f)}$. This gives us a way to estimate the size of the Julia sets of class $\b$ functions from below by considering their escaping sets.

\smallskip
Let now $f\in\b$. By definition, we can choose $R_0>0$ large enough such that $\sing(f^{-1})\subset D(0,R_0)$. Put $B(R_0):=\hc\setminus\overline{D(0,R_0)}$. Then a result of Rippon and Stallard \cite[Lemma 2.1]{rippon6} says that every component of $f^{-1}(B(R_0))$ is simply connected. For entire $f$ this is due to Eremenko and Lyubich \cite{eremenko2}. Moreover, if $\infty$ is not an asymptotic value, then every such component is bounded and contains exactly one pole of $f$.

\bigskip
\noindent\emph{Second order complex differential equations.}
To estimate the Hausdorff dimension of the escaping set for functions considered here, one crucial point is the asymptotic behaviour near poles. More precisely, we need to estimate the residues at poles. For Nevanlinna functions, this will be deduced from the connection between these functions and the second order complex differential equations, which is known for quite a long time. We refer to \cite[Chapter 5]{hille3} or \cite[Chapter 4]{langley1} for details.

\medskip
Let $f$ be a Nevanlinna function. Without loss of generality, we assume that the Schwarzian derivative of $f$ is
\begin{equation}\label{sd}
\s_f(z)=2p(z),
\end{equation}
where $p$ is a polynomial of degree $m$. The relation between \eqref{sd} and the following second order complex differential equation
\begin{equation}\label{so}
w''(z)+p(z)w(z)=0
\end{equation}
is as below.  If $w_1$ and $w_2$ are two linearly independent solutions of \eqref{so} then the Wronskian $W(w_1,w_2)$ is a non-zero constant. Moreover, the ratio $f:=w_2/w_1$ satisfies the above Schwarzian differential equation \eqref{sd}.

By making use of the Liouville transformation
$$W(Z)=p(z)^{1/4}w(z),\,\,\, Z=\int^{z}p(s)^{1/2}ds,$$
one can transfer \eqref{so} to a sine-type equation
\begin{equation}\label{so1}
W''(Z)+\left(1-F(Z)\right)W(Z)=0,
\end{equation}
where
$$F(Z)=\frac{1}{4}\,\frac{p''(z)}{p(z)^2}-\frac{5}{16}\,\frac{p'(z)^2}{p(z)^3}.$$
What is important here is that $F(Z)=\mathcal{O}(1)/Z^{2}$ for $|Z|$ large and for $Z$ in the following sector
\begin{equation}\label{psec}
\left\{Z:~|\arg Z|<\pi-\delta~\right\},
\end{equation}
where $0<\delta<\pi$ is some constant. There are two linearly independent solutions $W_1$ and $W_2$ of \eqref{so1} which have the following asymptotics
$$W_1(Z)= e^{iZ}(1+\varepsilon_1(Z)),~~\,~W_2(Z)= e^{-iZ}(1+\varepsilon_2(Z))$$
and $|\varepsilon_j(Z)|=\mathcal{O}(1/|Z|)$ for $Z$ in the sector \eqref{psec}. This gives two linearly independent solutions
$$w_1:=p(z)^{-1/4}W_1~\,\,\text{and}\,\,~w_2:=p(z)^{-1/4}W_2$$
of the equation \eqref{so} in a sector
$$S:=\left\{z: \left|\arg z-\theta_0\right|<\frac{2\pi}{m+2}-\delta'~\text{and~~$|z|$~large}\right\}$$
where $0<\delta'<\frac{2\pi}{m+2}$ is a constant depending on $\delta$. Here $\theta_0$ is a solution of
\begin{equation}\label{cr}
\arg a_m +(m+2)\theta=0~\,(\mod 2\pi),
\end{equation}
where $a_m$ is the leading coefficient of $p(z)$. Every ray with argument $\theta$ solving the above equation \eqref{cr} is called a critical ray or a Julia ray. The solutions $w_1$ and $w_2$ have no zeros in the sector $S$. But it is known that, if $\lambda,\mu\neq 0$ are constants, then a solution of the form $\lambda w_1+ \mu w_2$ will have zeros near the critical ray.

Thus, for suitable constants $a,\,b,\,c,\,d\in\c$ with $ad-bc\neq 0$ the function $f$ can be represented as
\begin{equation*}
f(z)=\frac{a\,w_1(z)+b\,w_2(z)}{c\,w_1(z)+d\,w_2(z)}=\frac{a\,W_1(Z)+b\,W_2(Z)}{c\,W_1(Z)+d\,W_2(Z)}=:\frac{a\,g(Z)+b}{c\,g(Z)+d},
\end{equation*}
where $g(Z)=W_1(Z)/W_2(Z)\sim e^{i2Z}$ for large $Z$ in the sector \eqref{psec}. Here and in the following $A\sim B$ means that $A$ is equal to $B$ up to some constant. Thus $\log g(Z)\sim i2Z$ and hence
$$\frac{g'(Z)}{g(Z)}\sim 2i$$
for large $Z$ in the sector \eqref{psec}.
Now let $z_0$ be a pole of $f$, that is, $z_0$ is a zero of $cg(Z)+d$. We will compute the residue of $f$ at $z_0$. Suppose that $Z_0=\int^{z_0}p(s)^{1/2}ds$. This gives $cg(Z_0)+d=0$. Then
\begin{equation*}
\begin{aligned}
\res(f,z_0)&=\lim_{z\to z_0}(z-z_0)f(z)\\
&=\lim_{z\to z_0}(z-z_0)\,\frac{a\,g(Z)+b}{c\,g(Z)+d}\\
&=(a\,g(Z_0)+b)\lim_{z\to z_0}\frac{(z-z_0)}{c\,g(Z)+d}\\
&=(a\,g(Z_0)+b)\,\frac{1}{c\,g'(Z_0)\,p(z_0)^{1/2}}\\
&\sim\left(-a\,\frac{d}{c}+b\right)\frac{1}{c\cdot 2ig(Z_0)\,p(z_0)^{1/2}}\\
&=\frac{1}{2i}\left(\frac{a}{c}-\frac{b}{d}\right)\frac{1}{p(z_0)^{1/2}}\\
&\sim c'\cdot z_{0}^{-m/2},
\end{aligned}
\end{equation*}
where $c'$ is some constant.

\bigskip
\noindent\emph{Nevanlinna functions.} We collect here some properties of Nevanlinna functions for future purposes. Let $f$ be such a function. Then $f$ has only finitely many asymptotic values and only simple poles (if there are any). The order $\rho(f)$ of $f$ is $(m+2)/2$, where $m$ is the degree of the Schwarzian derivative of $f$, that is, the degree of the polynomial $p$ in \eqref{sd}, see, for instance, \cite[Chapter XI, Section 3]{nevanlinna6}. Suppose that $a_j$ are poles of $f$, ordered in such a way that $\dots\leq |a_j|\leq |a_{j+1}|\leq\dots$. Then we have an asymptotic formula for $|a_j|$ for sufficiently large $j$. That is,
$$|a_{j}|\sim c'' j^{1/\rho(f)}=c'' j^{2/(m+2)},$$
where $c''$ is a constant, see \cite[Chapter 4]{langley1} for instance.

\medskip
Suppose from now that $f$ is a Nevanlinna function for which $\infty$ is not an asymptotic value. As noted at the beginning of this section, for large $R>R_0$ every component $U$ of $f^{-1}(B(R))$ is a bounded, simply connected domain containing exactly one pole of $f$. Since $f$ has only simple poles, the restriction of $f$ on each $U$ is actually a conformal map onto $B(R)$. The component containing the pole $a_j$ will be denoted by $U_j$ and the inverse branch of $f$ into $U_j$ is denoted by $g_j$. For simplicity, we set $b_j=\res(f,a_j)$. Thus, the above calculation of residues at poles yields that
\begin{equation}\label{residue}
|b_j|\sim c'\cdot |a_{j}|^{-m/2}\sim c\cdot j^{-m/(m+2)}.
\end{equation}

We will estimate the (Euclidean) derivative of $g_j$ and the size of $U_j$. Since $U_j$ is simply connected, there is a conformal map $\varphi_j : U_j \to D(0, R^{-1})$ such that $\varphi_{j}(a_j)=0$ and $\varphi'_{j}(a_j)=1/b_j$. It follows from the Koebe's one-quarter theorem (\cite[Theorem 1.3]{pommerenke1}) that
$$U_j=\varphi_{j}^{-1}\left(D\left(0,R^{-1}\right)\right)\supset D\left(a_j, \frac{|b_j|}{4R}\right).$$
Moreover, by choosing $R>0$ appropriately one may, by using Koebe's distortion theorem again, achieve that
$$U_j\subset D\left(a_j,  \frac{2|b_j|}{R} \right).$$
For details, we refer to \cite[Section 3]{bergweiler2} or \cite[Proof of Theorem 3.1]{cuiwei1}. Therefore, altogether we have the following estimate of the size of $U_j$ for large $j$:
\begin{equation}\label{distortion}
D\left(a_j, \frac{|b_j|}{4R}\right)\subset U_j\subset D\left(a_j,  \frac{2|b_j|}{R} \right).
\end{equation}
This means that
\begin{equation}\label{diameter}
\diam U_j \leq\frac{4|b_j|}{R}.
\end{equation}
In particular, if $j$ is large enough, then $U_j$ will be contained in $B(R)$. Moreover, since
$$g_j(z)=\varphi_{j}^{-1}\left(\frac{1}{z} \right)$$
we see that
$$|g'_{j}(z)|\leq \frac{12\left|\left(\varphi_{j}^{-1}\right)'(0)\right|}{|z|^2}=\frac{12|b_j|}{|z|^2}$$
for $z$ in any simply connected domain $D\subset B(R)\setminus\{\infty\}$. In particular, if we choose the domain $D$ to be $U_k$ for large $k$, then we can have
\begin{equation}\label{derivative}
|g'_{j}(z)|\leq \frac{12|b_j|}{|z|^2}\leq C\cdot \frac{|b_j|}{|a_k|^2}
\end{equation}
for some constant $C>0$ independent of $j$ and $k$.

\section{The estimate of the Hausdorff dimension}

In this section, we will compute the Hausdorff dimension of the escaping set from above and from below separately.

\bigskip
\noindent\emph{The upper bound.} Now for $U_j\subset B(R)$, using \eqref{distortion} and \eqref{diameter} we can estimate the diameter of $g_{j}(U_k)$ as follows.
\begin{equation*}
\begin{aligned}
\diam g_{j}(U_k) &\leq \sup_{z\in U_k}|g'_{j}(z)|\diam U_k\\
&\leq C\cdot \frac{|b_j|}{|a_k|^2}\cdot \frac{4|b_k|}{R}.
\end{aligned}
\end{equation*}
Here $C$ is some constant. Inductively one can estimate the successive pullbacks of $U_k$ under finitely many inverse branches of $f$. More precisely, suppose that $U_{j_1},\, U_{j_2},\,\dots\, U_{j_l} \subset B(R)$. We have
\begin{equation}\label{diad}
\begin{aligned}
\diam \left(g_{j_1} \circ g_{j_2} \circ \dots \circ g_{j_{l-1}}\right)(U_{j_l}) &\leq C\frac{|b_{j_1}|}{|a_{j_2}|^2}\cdots C\frac{|b_{j_{l-2}}|}{|a_{j_{l-1}}|^2} \cdot C\frac{|b_{j_{l-1}}|}{|a_{j_{l}}|^2} \cdot \frac{4|b_{j_l}|}{R}\\
&=C^{l-1}\,\frac{4}{R}\cdot|b_{j_1}|\cdot\frac{|b_{j_2}|}{|a_{j_2}|^2}\cdots \frac{|b_{j_{l-1}}|}{|a_{j_{l-1}}|^2}\cdot\frac{|b_{j_l}|}{|a_{j_l}|^2}\\
&=C^{l-1}\,\frac{4}{R}|b_{j_1}|\,\prod_{k=2}^{l}\,\frac{|b_{j_k}|}{|a_{j_k}|^2}.
\end{aligned}
\end{equation}
For any set $K\subset U_k$, the relation between Euclidean and spherical diameters of $K$ is as follows:
$$\diam_{\chi}(K)\leq \frac{C_1}{|a_k|^2}\diam (K),$$
where $C_1$ is some constant. We omit details here. Thus, we have
\begin{equation}\label{diame}
\diam_{\chi} \left(g_{j_1} \circ g_{j_2} \circ \dots \circ g_{j_{l-1}}\right)(U_{j_l}) \leq C^{l-1}\,\frac{4C_1}{R}\,\prod_{k=1}^{l}\,\frac{|b_{j_k}|}{|a_{j_k}|^2}.
\end{equation}

To find an upper bound for the dimension of the escaping set, we need to construct a cover of $\I(f)$. For this purpose, we define
$$\j_{R}(f):=\left\{z\in B(R): f^{n}(z)\in B(R)~\text{for all}~n\in\mathbb{N}\right\}$$
which is exactly the set we defined in the introduction, and
$$\I_{R}(f)=\j_{R}(f)\cap\I(f).$$
Suppose that $E_{l}$ is the collection of all components $V$ of $f^{-l}(B(R))$ such that $f^{k}(V)\subset B(R)$ for $0\leq k\leq l-1$. It is not difficult to see that $E_{l}$ is a cover of 
$$\left\{z\in B(3R): f^{n}(z)\in B(3R)~\text{for}~0\leq k\leq l-1\right\}.$$
Now if we take $t>(m+2)/(m+4)$, then, by using \eqref{residue}
\begin{equation*}
\begin{aligned}
\sum_{V\in E_l}(\diam_{\chi} V)^t &\leq \left(C^{l-1}\,\frac{4C_1}{R}\right)^t \sum_{j_1=M}^{\infty}\cdots\sum_{j_l=M}^{\infty}\prod_{k=1}^{l}\left(\frac{|b_{j_k}|}{|a_{j_k}|^2}\right)^t\\
&\leq \left(C^{l-1}\,\frac{4C_1}{R}\right)^t \,\prod_{k=1}^{l}\left( \sum_{j_k=M}^{\infty}j_{k}^{-(m+4)/(m+2)}\right)^t\\
&<\infty,
\end{aligned}
\end{equation*}
for $R$ large, where $M$ is the largest number $j$ such that $U_j\not\subset B(R)$. Thus $\dim\j_{3R}(f)\leq t$, which implies that $\dim\I_{3R}(f)\leq t$. Therefore, we have
$$\dim\I(f)\leq \frac{m+2}{m+4}.$$

\bigskip

\bigskip

\noindent\emph{The lower bound.} The estimate of the lower bound of the escaping set will use a result of McMullen \cite{mcmullen11}. To describe his result, assume that $E_{\ell}$ is a collection of disjoint compact subsets of $\hc$ for each $\ell\in\N$ satisfying
\begin{itemize}
\item[(i)] each element of $E_{\ell+1}$ is contained in a unique element of $E_{\ell}$;
\item[(ii)] each element of $E_{\ell}$ contains at leat one element of $E_{\ell+1}$.
\end{itemize}
Suppose that $\overline{E}_{\ell}$ is the union of all elements of $E_{\ell}$. Let $E=\bigcap \overline{E}_{\ell}$. Suppose that for $V\in E_{\ell}$,
$$\dens_{\chi}(\overline{E}_{\ell+1}, V)\geq \Delta_{\ell}$$
and
$$\diam_{\chi} V \leq d_\ell$$
for two sequences of positive real numbers $(\Delta_{\ell})$ and $(d_{\ell})$. Then McMullen's result says that
$$\dim E \geq 2-\limsup_{\ell\to\infty}\frac{\sum_{j=1}^{\ell+1}|\log\Delta_j|}{|\log d_{\ell}|}.$$

\medskip
Now we start to estimate the lower bound. Basically we will follow the approach of Bergweiler and Kotus in \cite{bergweiler2}. In the following, $B, B_1, B_2,\dots$ will denote constants.

Suppose that $E_{\ell}$ is the same as in the estimate of the upper bound. Then $E\subset \j_{R}(f)$. So if $V\in E_{\ell}$ is such that $f^{k}(V)\subset U_{j_{k+1}}$ for $0\leq k\leq {\ell}-1$, then it follows from \eqref{diad} that
\begin{equation}\label{die}
\diam_{\chi}(V)\leq C^{\ell-1}\,\frac{4C_1}{R}\,\prod_{k=1}^{\ell}\,\frac{|b_{j_k}|}{|a_{j_k}|^2}\leq \left(\frac{B_1}{R^{m/2+2}} \right)^{\ell}=:d_{\ell}
\end{equation}
for certain constant $B_1$.

On the other hand, put, for $s\geq R$ large,
$$A(s):=\{z: s<|z|<2s\}.$$
We first estimate $\dens(\overline{E}_1, A(s))$. Note that $\area A(s)= 3\pi s^2$. The number of $U_j$'s contained in $A(s)$ is  greater than or equal to $B_2\cdot \dfrac{s}{\max_j \diam U_j}$, $\max$ ranges over all $j$ such that $U_j$ is contained in $A(s)$. Thus it follows from \eqref{residue} and \eqref{distortion} that this number is no less than $B_3\,s^{m/2+1}$.

For those $U_j$ which are contained in $A(s)$, it follows from \eqref{distortion} and \eqref{residue} that
$$\diam U_j\geq\, B_4\, \frac{|b_j|}{2R}\,\geq\, B_5\, \frac{|a_j|^{-m/2}}{2R}\,\geq\frac{B_6}{R}\cdot s^{-m/2}\geq B_6\, \frac{1}{s^{m/2+1}}.$$
Thus we have
$$\area \left(\overline{E}_1 \cap A(s)\right)\geq B_3 s^{m/2+1}\cdot\pi \left(\frac{B_6}{2s^{m/2+1}}\right)^2 =\frac{\pi B_{6}^{2}B_3}{4}\frac{1}{s^{m/2+1}}=:\frac{B_7}{s^{m/2+1}},$$
which yields that
$$\dens\left(\overline{E}_1, A(s)\right)\geq\frac{B_7}{3\pi}\frac{1}{s^{m/2+3}}.$$

Since $V\in E_{\ell}$, by definition $f^{\ell-1}(V)=U_j$ for some $j$, where $U_j$ is a component of $f^{-1}(B(R))$ contained in $B(R)$. Thus $f^{\ell-1}(\overline{E}_{\ell}\cap V)=\overline{E}_{2}\cap U_j$. Thus, similarly as in \cite[pp. 5384]{bergweiler2}, by using Koebe's distortion theorem again, we have
$$\dens\left(\overline{E}_{\ell+1}, V\right)\geq B_8\cdot \dens \left(\overline{E}_{2}, U_j\right)\geq B_9\cdot \dens \left(\overline{E}_{1}, A(s)\right)\geq \frac{B_7 B_9}{3\pi}\frac{1}{s^{m/2+3}}.$$

Take $s=2^k R$, using \cite[(4.2)]{bergweiler2} we have that
\begin{equation}\label{dene}
\dens_{\chi}\left(\overline{E}_{\ell+1}, V\right)\geq\frac{B}{R^{m/2+3}}=:\Delta_{\ell},
\end{equation}
where $B>0$ is some constant. Thus, with \eqref{die} and \eqref{dene} and using McMullen's result mentioned above, we have
$$\dim E\geq 2-\limsup_{\ell\to\infty}\frac{(\ell+1)\left(\log B-(m/2+3)\log R\right)}{\ell\,\left(\log B_1-(m/2+2)\log R\right)}=2-\frac{\log B-(m/2+3)\log R}{\log B_1-(m/2+2)\log R}.$$
Then by taking $R\to\infty$, we obtain that
$$\dim E \geq \frac{m+2}{m+4}.$$
Thus we have that 
$$\dim\j_{R}(f)\geq \frac{m+2}{m+4}.$$
Combined with previous estimate on the upper bound we just proved that
$$\lim_{R\to\infty}\dim\j_{R}(f)=\frac{m+2}{m+4}.$$

For the estimate of the lower bound for the Hausdorff dimension of the escaping set, instead of taking a fixed $R$, we now consider an increasing sequence $(R_l)$ tending to $\infty$. Then for each $l$ we define $E_l$ as the collection of components of $f^{-l}(B(R))$ satisfying $f^{k}(V)\subset B(R_l)$ for $0\leq k\leq l-1$. Denote by $\overline{E}_l$ the union of the components in $E_l$ and set $E=\cap_{l}\overline{E}_l$. Then it follows that $E\subset \I(f)$. Similarly as above computations one can obtain that
$$\dim E\geq \frac{m+2}{m+4}.$$
Thus we have
$$\dim\I(f)\geq \frac{m+2}{m+4}.$$

\bigskip
In conclusion we have proved that the Hausdorff dimension of the escaping set of a Nevanlinna functions is equal to ${(m+2)}/{(m+4)}$.

\begin{ack}
I would like to thank Walter Bergweiler, Weiyuan Qiu, Weixiao Shen and Jun Wang for useful comments and valuable discussions. I am also grateful for the helpful comments of the referees.
\end{ack}

\bigskip
\emph{Shanghai Center for Mathematical Sciences,  Fudan University,  2005 Songhu Road, Shanghai 200438,  China}

\medskip
\emph{cuiweiwei@fudan.edu.cn}

\end{document}